# On Algebraic Solutions of Polynomial Equations of Degree n in one Variable.

Gerry Martens
Math Research GH Consulting Fuerstenfeldbruck Germany
GerryMrt@aol.com

**Abstract**

We will show that the roots of a polynomial equation in one variable of degree $n$ are related to the solutions of a symmetric quadratic form in $n-1$ variables with constant positive integer coefficients. The classic polynomial notation will be rewritten to define a characteristic discriminant of a polynomial of degree $n$. A new set of characteristic roots allows expressing the characteristic discriminant as the result of a symmetric quadratic form.

## 1) Rewriting the Classic Polynomial

The polynomial of a degree $n$ in one variable

(1) $$f(x) = a_0 + a_1 x + .. + a_{n-2} x^{n-2} + a_{n-1} x^{n-1} + a_n x^n = \sum_{i=0}^{n} a_i x^i$$

with coefficients $a_i$ being integer, real, or complex can be rewritten as

(2) $$f(x) = \frac{(n! a_n x + (n-1)! a_{n-1})^n - \sum_{i=2}^{n} \binom{n}{i} ((n-1)! a_{n-1})^i (n! a_n)^{n-i} x^{n-i} + n!^n a_n^{n-1} \sum_{i=0}^{n-2} a_i x^i}{n!^n a_n^{n-1}}$$

by completing the $n$-th power of the $(n-1)$-th derivative of the the classic notation(1).

This can be seen as

$$f(x) = \frac{(f^{(n-1)}(x))^n - g(x) + n!^n a_n^{n-1} \sum_{i=0}^{n-2} a_i x^i}{n!^n a_n^{n-1}}$$

with $(f^{(n-1)}(x))^n = (n! a_n x + (n-1)! a_{n-1})^n$ the $n$-th power of the $(n-1)$-th derivative and $g(x) = \sum_{i=2}^{n} \binom{n}{i} ((n-1)! a_{n-1})^i (n! a_n)^{n-i} x^{n-i}$ the terms of the derivative function

$(f^{(n-1)}(x))^n$ not part of the polynomial (Here $\binom{n}{i}$ is the binomial coefficient)

and $n!^n a_n^{n-1} \sum_{i=0}^{n-2} a_i x^i$ the polynomial terms of powers $(n-2)$ down to 0 multiplied by

$n!^n a_n^{n-1}$. Canceling the common factor $n!^n a_n^{n-1}$ we get the classic polynomial.

Gerry Martens Ferdinand Miller Strasse 36 D-82256 Fuerstenfeldbruck Germany +49 160 94767547

On Algebraic Solutions of Polynomial Equations of Degree n in one Variable.

## 2) The Polynomial Equation

The polynomial equation can be expressed as
(3) $$f(x) = a_n(x+x_1)(x+x_2)(x+x_3)...(x+x_i)...(x+x_n) = 0$$
Using representation (2) we can rewrite the polynomial equation (3) showing one part left and one part right of the equal sign. With the common factor cancelled out we get
(4) $$(n!a_n x + (n-1)!a_{n-1})^n = \sum_{i=2}^{n} \binom{n}{i}((n-1)!a_{n-1})^i (n!a_n)^{n-i} x^{n-i} - n!^n a_n^{n-1} \sum_{i=0}^{n-2} a_i x^i$$
According to the fundamental theorem of algebra which also applies to (4) both sides will only be equal for the roots $-x_1, -x_2, -x_3, ... -x_i, ... -x_n$ of the polynomial.

## 3) The Characteristic Discriminant of the Polynomial

By differentiation we can generate $(n-2)$ derivative equations from equation (4) (The quadratic equation is a special case and can be used without differentiation). Each of these equations will have a perfect power of descending degree $(n-1)$ down to 2 on the left side. The characteristic discriminant equation is the $(n-2)$-th derivative of the polynomial equation (4) which we can express for the

Quadratic equation as: $n=2$   $(2a_n x + a_{n-1})^2 = a_{n-1}^2 - 4a_n a_{n-2}$
Cubic equation as:   $n=3$   $(6a_n x + 2a_{n-1})^2 = 4a_{n-1}^2 - 12a_n a_{n-2}$
Quartic equation as:   $n=4$   $(24a_n x + 6a_{n-1})^2 = 36a_{n-1}^2 - 96a_n a_{n-2}$
Quintic equation as:   $n=5$   $(120a_n x + 24a_{n-1})^2 = 576a_{n-1}^2 - 1440a_n a_{n-2}$
Sextic equation as:   $n=6$   $(720a_n x + 120a_{n-1})^2 = 14400a_{n-1}^2 - 34560a_n a_{n-2}$
Septic equation as:   $n=7$   $(5040a_n x + 720a_{n-1})^2 = 518400a_{n-1}^2 - 1209600a_n a_{n-2}$
Octic equation as:   $n=8$   $(40320a_n x + 5040a_{n-1})^2 = 25401600a_{n-1}^2 - 58060800a_n a_{n-2}$

And the general characteristic discriminant equation for a polynomial of degree n can be expressed as:
(5) $$(n!a_n x + (n-1)!a_{n-1})^2 = (n-1)!^2 a_{n-1}^2 - 2n!(n-2)!a_n a_{n-2}$$

The right part of this equation is what we call the characteristic discriminant $D_n$ of the polynomial of degree $n$.
(6) $$D_n = (n-1)!^2 a_{n-1}^2 - 2n!(n-2)!a_n a_{n-2}$$

As we can see $D_n$ only depends on the coefficients $a_n, a_{n-1}, a_{n-2}$ of the highest degree terms and has a divisor $(n-1)!(n-2)!$. The coefficients $a_n, a_{n-1}, a_{n-2}$ will determine whether the value is positive negative, square, zero, integer, real or complex.

Solving equation (5) for $x$ one finds two solutions, which are the roots in case of the Quadratic polynomial. For the higher degree polynomials one of these solutions only represents a root if the polynomial has $(n-1)$ equal roots.
The characteristic discriminant allows us to define a new set of $(n-1)$ roots which will be introduced next.

## 4) Representation of the Regular Roots, Reference- and Characteristic Roots

The factored form of the polynomial of degree $n$ as shown in (3) lists all the **regular roots**. Carrying out the multiplication we can derive the relationships between the regular roots and the coefficients $\dfrac{a_i}{a_n}$ of the polynomial $f(x) = \sum_{i=0}^{n} a_i x^i$.

The coefficient $a_{n-1}$ represents the sum of all regular roots when we assume $a_n = 1$. In this case we can write the coefficient $a_{n-1}$ for the following polynomial degrees as

$n = 2 \quad a_{n-1} = x_1 + x_2$
$n = 3 \quad a_{n-1} = x_1 + x_2 + x_3$
$n = 4 \quad a_{n-1} = x_1 + x_2 + x_3 + x_4$
$n = 5 \quad a_{n-1} = x_1 + x_2 + x_3 + x_4 + x_5$

For a polynomial of degree $n$ we can write

$$(7) \qquad a_{n-1} = (x_1 + x_2 + x_3 + \ldots x_i \ldots + x_{n-1} + x_n) = \sum_{i=1}^{n} x_i$$

A new set of $(n-1)$ variables $b_1$ to $b_{n-1}$ which will be referenced as the **characteristic roots** are defined for the following polynomial degrees as

$n = 2 \quad x_2 = x_1 + b_1$
$n = 3 \quad x_2 = x_1 + b_1, \quad x_3 = x_2 + b_1 + b_2$
$n = 4 \quad x_2 = x_1 + b_1, \quad x_3 = x_2 + b_1 + b_2, \quad x_4 = x_3 + b_1 + b_2 + b_3$
$n = 5 \quad x_2 = x_1 + b_1, \quad x_3 = x_2 + b_1 + b_2, \quad x_4 = x_3 + b_1 + b_2 + b_3, \quad x_5 = x_4 + b_1 + b_2 + b_3 + b_4$

For a polynomial of degree $n$ we can write a set of $n-1$ equations as

$$(8) \qquad x_i = x_{i-1} + \sum_{j=1}^{i-1} b_j \quad \text{for } i = 2 \text{ to } n.$$

Substituting the higher roots $x_{i>1}$ in (8) and combining this result with our first equation (7) for the coefficient $a_{n-1}$ we obtain an equation for $x_1$

$$(9) \qquad x_1 = \frac{1}{n}(a_{n-1} - \sum_{i=1}^{n-1} \frac{1}{2}((n-i)^2 + (n-i))b_i).$$

The set of equations (8) and equation (9) together will allow us to yield the regular roots $x_i$ once we have found a solution for our characteristic roots $b_i$. The polynomial can this way be represented by one regular root say $x_1$ defined as the **reference root** and the characteristic roots $b_1, b_2, b_3, \ldots b_{n-1}$.

Looking at the regular roots as an ordered list we can define alternate sets of characteristic roots by rearranging the regular roots. The number of permutations on this list is $n!$ defining the number of sets of characteristic roots.
This new definition allows expressing the polynomial coefficients as functions of the reference root and the characteristic roots.
Geometrically the reference root is a slider for a polynomial with constant characteristic roots.

## 5) The Algebraic Notation of the Characteristic Roots

The characteristic discriminant (6) will allow us to construct a symmetric quadratic form of the characteristic roots with constant positive integer coefficients for a polynomial of degree $n$. These constant positive integer coefficients are not function of the polynomial coefficients or regular roots. The quadratic polynomial is a special case only having one variable and coefficient 1. To obtain these symmetric quadratic forms we first substitute the coefficients $a_n, a_{n-1}, a_{n-2}$ in (6) using the root relationship equations for the regular roots $x_i$ and secondly we substitute the regular roots using equations (8) for the characteristic roots $b_i$. Combining terms we get for a polynomial of degree

$n = 2$
$$D_n = a_{n-1}^2 - 4a_n a_{n-2} = (a_n(x_1 + x_2))^2 - 4a_n(x_1 x_2) = a_n^2(x_1 - x_2)^2$$
$$D_n = a_n^2 b_1^2$$

$n = 3$
$$D_n = 4a_{n-1}^2 - 12 a_n a_{n-2} = 4(a_n(x_1 + x_2 + x_3))^2 - 12 a_n(a_n(x_1(x_2 + x_3) + x_2 x_3))$$
$$D_n = 2a_n^2(6b_1^2 + 2b_2^2 + 2.3 b_1 b_2)$$

$n = 4$
$$D_n = 36 a_{n-1}^2 - 96 a_n a_{n-2} = 36(a_n(x_1 + x_2 + x_3 + x_4))^2 - 96 a_n(a_n(x_1(x_2 + x_3 + x_4) + x_2(x_3 + x_4) + x_3 x_4)$$
$$D_n = 12 a_n^2 (20 b_1^2 + 11 b_2^2 + 3 b_3^2 + 2.14 b_1 b_2 + 2.6 b_1 b_3 + 2.5 b_2 b_3)$$

Developing this further for higher degrees we can conclude that the characteristic discriminant $D_n$ is equal to $(n-1)!(n-2)!$ multiplied by the square of the coefficient $a_n$ and a symmetric quadratic form with $n-1$ variables with constant positive integer coefficients.
The symmetric quadratic form has a matrix notation

(10) $$B_n^T H_n B_n = \sum_{i=1}^{n-1} \sum_{j=1}^{n-1} h_{ij} b_i b_j$$

Where $B_n^T$ is the transposed column vector of the characteristic roots $b_1, b_2, b_3, \ldots b_{n-1}$,

$H_n$ represents a Hermitian (symmetric) matrix of $(n-1)*(n-1)$ elements with

$h_{ij} = h_{ji}$ being the coefficients which are constant for the polynomial degree $n$,

$B_n$ is the normal column vector of the characteristic roots $b_1, b_2, b_3, \ldots b_{n-1}$.

Rewriting the quadratic form to the matrix notation we get

$n = 2 \qquad D_2 = (n-1)!(n-2)! a_2^2 B_2^T H_2 B_2 \qquad H_2 = [1]$

$n = 3 \qquad D_3 = (n-1)!(n-2)! a_3^2 B_3^T H_3 B_3 \qquad H_3 = \begin{bmatrix} 6 & 3 \\ 3 & 2 \end{bmatrix}$

$n = 4 \qquad D_4 = (n-1)!(n-2)! a_4^2 B_4^T H_4 B_4 \qquad H_4 = \begin{bmatrix} 20 & 14 & 6 \\ 14 & 11 & 5 \\ 6 & 5 & 3 \end{bmatrix}$

We notice that the characteristic discriminant can be represented by the coefficients $a_n, a_{n-1}, a_{n-2}$ which are function of $x_1, x_2, x_3, \ldots x_n$ or the characteristic roots $b_1, b_2, b_3, \ldots b_{n-1}$.

## 6) The Matrix Notation of the Characteristic Discriminant.

For a polynomial of degree $n$ the matrix notation of the characteristic discriminant $D_n$ can be written as

(11) $\quad D_n = (n-1)!(n-2)!a_n^2 B_n^T H_n B_n = (n-1)!^2\, a_{n-1}^2 - 2n!(n-2)!a_n a_{n-2}$

This equation shows that $D_n$ has a divisor $(n-1)!(n-2)!a_n^2$

Dividing $D_n$ by $(n-1)!(n-2)!a_n^2$ we get

(12) $\quad \dfrac{D_n}{(n-1)!(n-2)!a_n^2} = B_n^T H_n B_n = \dfrac{1}{a_n^2}((n-1)a_{n-1}^2 - 2na_n a_{n-2})$

The general structure of the matrix $H_n$ can be written as

$$H_n = \begin{bmatrix} \frac{1}{12}(n^4 - n^2) & & & & & & & & & h_{1(n-1)} \\ & h_{22} & & & & & & & & \ldots \\ & & h_{33} & & & & & & & \ldots \\ & & & \ldots & & & & & & \ldots \\ & & & & h_{ij} & & & & & 6n-21 \\ & & & & & \ldots & & & & 5n-15 \\ & & & & & & \ldots & & & 4n-10 \\ & & & & & & & \ldots & & 3n-6 \\ & & & & & & & & \ldots & 2n-3 \\ h_{(n-1)1} & \ldots & \ldots & \ldots & 6n-21 & 5n-15 & 4n-10 & 3n-6 & 2n-3 & n-1 \end{bmatrix}$$

The elements of $H_n$ can be generated by the following lower triangular matrix.

(13) $\quad h_{ij} = \begin{cases} h_{ij} & i \geq j \\ 0 & i < j \end{cases}$

$\quad h_{ij} = \dfrac{1}{6}(ni^3 - \dfrac{3}{2}n(n+1)i^2 + \dfrac{1}{2}n(3n+1)i + \dfrac{1}{2}(n^4 - n^2)) - \dfrac{1}{4}(i^2 - (2n+1)i + n(n+1))(j-1)j$

To complete the matrix $H_n$ we copy the elements from the lower triangular matrix $h_{ij} = h_{ji}$ over to the upper part.
The elements $h_{ij}$ of the principal diagonal for which $i = j$ can be expressed as the following Quartic polynomial

(14) $\quad h_{ii} = \dfrac{1}{12}(-3i^4 + 2(4n+3)i^3 - 3(2n^2 + 4n + 1)i^2 + 2n(3n+2)i + n^4 - n^2)$

$\quad\quad\quad$ for $i = 1$ to $(n-1)$

This Quartic has the following factored form

(15) $\quad h_{ii} = -\dfrac{1}{4}(i - n - 1)(i - n)(i - \dfrac{1}{6}(2n + 3 + \sqrt{16n^2 + 9}))(i - \dfrac{1}{6}(2n + 3 - \sqrt{16n^2 + 9}))$

## 7) The Characteristic Discriminant for Polynomials up to Degree 8.

Below the characteristic discriminant notation for the first eight degrees:

$$D_2 = a_2^2 B_2^T H_2 B_2 = a_1^2 - 4a_2 a_0 \qquad H_2 = [1]$$

$$D_3 = 2a_3^2 B_3^T H_3 B_3 = 4a_2^2 - 12a_3 a_1 \qquad H_3 = \begin{bmatrix} 6 & 3 \\ 3 & 2 \end{bmatrix}$$

$$D_4 = 12 a_4^2 B_4^T H_4 B_4 = 36 a_3^2 - 96 a_4 a_2 \qquad H_4 = \begin{bmatrix} 20 & 14 & 6 \\ 14 & 11 & 5 \\ 6 & 5 & 3 \end{bmatrix}$$

$$D_5 = 144 a_5^2 B_5^T H_5 B_5 = 576 a_4^2 - 1440 a_5 a_3 \qquad H_5 = \begin{bmatrix} 50 & 40 & 25 & 10 \\ 40 & 34 & 22 & 9 \\ 25 & 22 & 16 & 7 \\ 10 & 9 & 7 & 4 \end{bmatrix}$$

$$D_6 = 2880 a_6^2 B_6^T H_6 B_6 = 14400 a_5^2 - 34560 a_6 a_4 \qquad H_6 = \begin{bmatrix} 105 & 90 & 66 & 39 & 15 \\ 90 & 80 & 60 & 36 & 14 \\ 66 & 60 & 48 & 30 & 12 \\ 39 & 36 & 30 & 21 & 9 \\ 15 & 14 & 12 & 9 & 5 \end{bmatrix}$$

$$D_7 = 86400 a_7^2 B_7^T H_7 B_7 = 518400 a_6^2 - 1209600 a_7 a_5$$

$$H_7 = \begin{bmatrix} 196 & 175 & 140 & 98 & 56 & 21 \\ 175 & 160 & 130 & 92 & 53 & 20 \\ 140 & 130 & 110 & 80 & 47 & 18 \\ 98 & 92 & 80 & 62 & 38 & 15 \\ 56 & 53 & 47 & 38 & 26 & 11 \\ 21 & 20 & 18 & 15 & 11 & 6 \end{bmatrix}$$

$$D_8 = 3628800 a_8^2 B_8^T H_8 B_8 = 25401600 a_7^2 - 58060800 a_8 a_6$$

$$H_8 = \begin{bmatrix} 336 & 308 & 260 & 200 & 136 & 76 & 28 \\ 308 & 287 & 245 & 190 & 130 & 73 & 27 \\ 260 & 245 & 215 & 170 & 118 & 67 & 25 \\ 200 & 190 & 170 & 140 & 100 & 58 & 22 \\ 136 & 130 & 118 & 100 & 76 & 46 & 18 \\ 76 & 73 & 67 & 58 & 46 & 31 & 13 \\ 28 & 27 & 25 & 22 & 18 & 13 & 7 \end{bmatrix}$$

## 8) Properties of the Characteristic Roots and Characteristic Discriminant

Using a variable $s$ as a second index for the roots we can address the permutations of the regular roots as $x_{is}$ and the characteristic roots as $b_{is}$.

A set of solutions is organized if the partial sums of the characteristic roots of a polynomial of degree $n \geq 2$ are all zero which we can write as

(16) $$\sum_{s=1}^{i!} b_{(i-1)s} = 0 \quad \text{for } i = 2 \text{ to } n$$

Regardless if the set is organized or not the sums of the characteristic roots over the complete set of solutions $n!$ are zero for every polynomial of degree $n \geq 2$. This we can write as

(17) $$\sum_{s=1}^{n!} b_{is} = 0 \quad \text{for } i = 1 \text{ to } n-1$$

The sums of the products of two characteristic roots of a set of solutions have the following relationship with the characteristic discriminant for every polynomial of degree $n > 2$.

(18) $$D_n = -\frac{1}{3}(n-1)! \sum_{s=1}^{n!} b_{1s} b_{2s}$$

$$D_n = (n-1)! \sum_{s=1}^{n!} b_{(1+i)s} b_{(3+i)s} \quad \text{for } i = 0 \text{ to } n-4 \text{ and } n \geq 4$$

$$D_n = -\frac{1}{4}(n-1)! \sum_{s=1}^{n!} b_{(2+i)s} b_{(3+i)s} \quad \text{for } i = 0 \text{ to } n-4 \text{ and } n \geq 4$$

All other sums of the products of two characteristic roots are zero

## 9) Conclusions

The fundamental theorem of algebra stating that

*"Every polynomial equation of degree n with complex coefficients has n roots in the complex numbers "*

can be stated in function of the characteristic roots as

*"Every polynomial equation of degree n with complex coefficients has n! sets of (n-1) characteristic roots in the complex numbers. The characteristic roots are the solutions of a symmetric quadratic form having (n-1) variables and constant positive integer coefficients for the polynomial degree n as described by the Hermitian matrix $H_n$ in section 6"*

This significant result allows us to conclude that finding the $n$ roots of a polynomial of degree $n$ is equivalent to finding one set of $n-1$ characteristic roots from a total of $n!$ sets.

## 10) Acknowledgements

The idea to write this paper came to me at the beginning of April looking through some old notes and books from my student time which my wife wanted to be moved because she wanted to bring a little lemon tree into my office at home. So many thanks to her.